\documentstyle[12pt]{article}

\begin{document}

\begin{center}
{\large \bf Ricci flow on compact K\"ahler manifolds of positive bisectional curvature}
\end{center}
\begin{center}
\begin{tabular}{ccc}
Huai-Dong Cao, Bing-Long Chen, and Xi-Ping Zhu\\
\end{tabular}
\end{center}

\vspace{0.4 cm}

\section{The main result}
We announce a new proof of the uniform estimate on the curvature of solutions to the Ricci flow on a compact K\"ahler manifold $M^n$ with positive bisectional curvature.\\

Given a compact K\"ahler manifold $M^n$ of complex dimension $n$ and a K\"ahler metric $\widetilde g=\sum \widetilde g_{i \bar j}(x)dz^i dz^{\bar j}$ on $M^n$ of positive holomorphic bisectional curvature, consider the normalized K\"ahler-Ricci flow
\begin{eqnarray}
{\partial \over {\partial t}}g_{i\bar j}(x,t)= -R_{i\bar j}(x,t) + g_{i\bar j}(x,t)
\end{eqnarray}
on $M^n$ with initial condition $g_{i\bar j}(x,0)=\widetilde g_{i\bar j}$. By a result of Bishop-Goldberg [2], we can assume the K\"ahler class of the metric $\widetilde g$ satisfies the condition:
\begin{eqnarray}
[\widetilde\omega]=[\widetilde\Sigma]=\pi c_1(M^n)
\end{eqnarray}
where $\widetilde\omega =({\sqrt{-1}/2}) \tilde g_{i \bar j}dz^i\wedge dz^{\bar j}$ and
$\widetilde\Sigma = ({\sqrt{-1}/2}) \tilde R_{i \bar j}dz^i\wedge dz^{\bar j}$
are the K\"ahler form, the Ricci form of the metric $\tilde g$ respectively, while
$c_1(M)$ denotes the first Chern class.
Under the normalized initial condition (2), the first author [3] (see also Proposition 1.1 in [4]) showed that the solution
$g(x, t)=\sum g_{i \bar j}(x,t)dz^i dz^{\bar j}$ to the normalized flow (1) exists for all time. Furthermore by the work of Mok [11] (and Bando [1] for $n=3$), the solution metric $g(x,t)$ is known to have positive bisectional curvature for any $t>0$. Our main result is
on the uniform estimate of the curvature independent of $t$:\\

\noindent {\bf Main Theorem} {\em Let $M^n$ be a compact K\"ahler manifold of complex dimension $n$ and $\widetilde g=\sum \widetilde g_{i \bar j}(x)dz^i dz^{\bar j}$ be a K\"ahler metric on $M^n$ of positive holomorphic bisectional curvature satisfying condition (2). Let $g(x,t)$ be the solution to the normalized K\"ahler-Ricci
flow (1) on $M^n$ with the initial condition $g(x,0)=\widetilde g$. Then $g(x,t)$ is nonsingular. Namely, the curvature of $g(x,t)$ is uniformly bounded independent of $t$. }\\

Whether the solution $g$ above is nonsingular or not has been one of the important open problems in the study of Hamilton's Ricci flow on compact K\"ahler manifolds. For $n=1$, the main theorem is proved by Richard Hamilton [9] based on his Harnack inequality and entropy estimate for the Gaussian curvature. (see also the further improvements
by B. Chow [8] and Hamilton [10]). In [6,7], Chen-Tian obtained the uniform estimate on the curvature of $g(x,t)$ based on a Moser-Trudinger type inequality in [16] for K\"ahler-Einstein manifolds and some new functionals.  Note the existence of K\"ahler-Einstein metrics on $M^n$ follows
from the solution to the Frankel conjecture by Mori [12] and Siu-Yau [15]. In contrast, our proof of the main theorem does not rely on the exsitence of K\"ahler-Einstein metrics, but instead on first author's Harnack inequality for the scalar curvature $R(x,t)$ which is a consequence of the Li-Yau-Hamilton estimate for the K\"ahler-Ricci flow [4], and a very recent local injectivity radius estimate of Perelman for the Ricci flow [13]. Thus, this gives a crucial estimate needed in a new proof of the Frankel conjecture via Hamilton's Ricci flow method. In next section, we'll outline the proof of the Main Theorem. Complete details of the proof and discussions about convergence of $g(x, t)$ as $t\rightarrow \infty$ will appear elsewhere.

\medskip
\noindent {\bf Remark} Richard Hamilton and the first author [5] are able to prove
the main theorem for $n=2$ by studying the geometry of singularity model from the blow up limits and certain dimension reduction argument.

\section{Outline of the proof}

Let $g(x,t)$ be the solution to the normalized K\"ahler-Ricci flow (1) on compact K\"ahler manifold $M^n$ with positive bisectional curvature whose initial metric satisfies condition (2).  It is easy to see that under flow (1), the volume $V=V(t)$, the total scalar curvature $\int_{M^n}R(x,t)d\nu_t$, and the average scalar curvature \[r=r(t)=\frac{1}{V}\int_{M^n}R(x,t)d\nu_t\] of $(M^n, g(x,t))$ are all constant in $t$. In fact, $r(t)=n$ for all $t$. Clearly, it suffices to show that the scalar curvature of $g$ is uniformly bounded from above in $t$:
 $$R(x,t)<C, \ \forall x\in M^n, \forall t\in [0,\infty)$$ for some constant $C>0$ indepedent of $t$. To do this, let us first recall \\

\noindent {\bf The Harnack inequlaity for the scalar curvature (Cao [4])} {\em For any
$x,y\in M^n$ and $0<t_1<t_2<\infty$, the scalar curvature $R(x,t)$ of solution metric $g(x,t)$ satisfies the inequlaity
\begin{eqnarray}
R(x,t_1)\leq R(y,t_2)\frac{e^{t_2}-1}{e^{t_1}-1}e^{{\Delta}/4}.
\end{eqnarray}
Here $\Delta = \Delta(x,y,t_1,t_2)=\inf_{\gamma}\int_{t_1}^{t_2} |\gamma'(\tau)|^2 d\tau$, where the infimum is taken over all space-time curves from $(x,t_1)$ to $(y, t_2)$, and $|\gamma'(\tau)|$ is the velocity of $\gamma$ at time $\tau \in [t_1, t_2]$.}

Now for any $x\in M^n$ and $t>1$, set $t_1=t$ and $t_2=t+1$. We can find $y\in M^n$ such that $R(y, t+1)=r(t+1)=n$. It then follows from the Harnack inequality (3)
\begin{eqnarray}
R(x,t)& \leq & R(y,t+1)\frac{e^{t+1}-1}{e^t-1}e^{\Delta(x,y,t,t+1)/4}\nonumber\\
 &\leq& n(e+1)e^{\Delta(x,y,t,t+1)/4}.
\end{eqnarray}
On the other hand, from (1) one can show that
\begin{eqnarray}
\Delta(x,y,t,t+1) \leq e^2 d^2(x,y;t)
\end{eqnarray}
where $d(x,y;t)$ denotes the geodesic distance between $x$ and $y$ with respect to the metric $g_{i\bar j}(t)$. In particular, by (4) and (5),
\begin{eqnarray*}
 R(z, t)\leq n(e+1)\exp(\frac 14 e^2)
\end{eqnarray*}
for all $z\in B_t(y,1)$, the geodesic ball centered at $y$ of radius $1$ with respect to
the metric $g_{i\bar j}(t)$. Then by Perelman's no local collapsing result (Theorem 4.1 or its corollary in [13]), there exists a constant $\beta>0$ independent of $t$ such that the volume of the geodesic ball $B_t(y,1)$ has a uniform lower bound:
\begin{eqnarray}
{\mbox Vol}(B_t(y, 1))\geq \beta.
\end{eqnarray}

From (6) one can, by a volume comparison argument of Yau (cf.
[14]), deduce a uniform upper bound on the diameter $d_t$ of
$(M^n, g(x,t))$: there is a constant $D=D(\beta)>0$ independent of
$t$, such that for all $t>1$
\begin{equation}
 d_t\leq D.
\end{equation}

But then, (4), (5), and (7) imply that, for any $t>1$ and any $x\in M^n$,
\begin{eqnarray*}
R(x,t)\leq n(e+1)\exp(e^2D^2/4),
\end{eqnarray*}
 which is our desired uniform estimate.
\medskip

Finally we should remark that the above argument also works for
compact K\"ahler manifolds with nonnegative bisectional curvature.

\noindent Institute for Pure and Applied mathematics at UCLA, and Texas A\&M University \\
{\em E-mail}:  hcao@ipam.ucla.edu or cao@math.tamu.edu
\medskip

\noindent Department of Mathematics, Zhongshang University, \& The
Institute of Mathematical Sciences, The Chinese University of Hong
Kong

\smallskip
\noindent {\em Email}: blchen@math.cuhk.edu.hk \ \& \ stszxp@zsu.edu.cn


\begin{thebibliography}{99}

\bibitem{Ba} S. Bando, {\em Compact $3$-folds with nonnegative bisectional curvature}, J. Differ. Geom.
{\bf 19} (1984), 283-297.

\bibitem{BG} R. L. Bishop and S. I. Goldberg, {\em On the second cohomologh group of a Kaehler manifold of positive curvature},  Proc. Amer. Math. Soc. {\bf 16} (1965), 119-122.


\bibitem{Ca85} H.-D. Cao, {\em Deformation of K\"ahler metrics to K\"ahler -Einstein metrics
on compact K\"ahler manifolds}, Invent. Math. {\bf 81} (1985), 359-372.


\bibitem{Ca92} H.-D. Cao, {\em On Hanarck's inequalities for the K\"ahler-Ricci
flow}, Invent. Math. {\bf 109} (1992), 247-263.

\bibitem{CH} H.-D. Cao and R. S. Hamilton, unpublished work.

\bibitem{CT1} X. Chen and G. Tian, {\em Ricci flow on K\"ahler-Einstein surfaces}, Invent. Math. {\bf 147} (2002), 487-544.

\bibitem{CT2} X. Chen and G. Tian, {\em Ricci flow on K\"ahler-Einstein manifolds}, preprint.

\bibitem{Ch}  B. Chow, {\em The Ricci flow on the $2$-sphere}, J. Differ. Geom. {\bf 33} (1991), 325-334.


\bibitem{Ha86} R. S. Hamilton,  {\em The Ricci flow on surfaces}, in Mathematics
 and General Relativity, Contemporary Mathematics {\bf 71} (1988), 237-261.

\bibitem{Hai} R. S. Hamilton, {\em An isoperimetric estimate for the Ricci flow on the 2-sphere}, Modern methods in complex analysis (Princeton, NJ, 1992), 201-222, Ann. of Math. Stud. {\bf 137}, Princeton Univ. Press, 1995.


\bibitem{Mo} N. Mok, {\em  The uniformization theorem for compact K\"ahler
manifolds of nonnegative holomorphic bisectional curvature}, J. Differ. Geom.
{\bf 27} (1988), 179-214.

\bibitem{Mr} S. Mori, {\em Projective manifolds with ample tangent bundles}, Ann. Math. {\bf 76} (1979), 213-234.

\bibitem{Pr} G. Perelman, {\em The entropy formula for the Ricci
flow and its geometric applications}, preprint.

\bibitem{ShY} R. Schoen and S.-T. Yau, {\em Lecture on Differential Geometry},
Conference Proceedings and Lecture Notes in Geometry and Topology, I. International Press, Cambridge, MA, 1994.


\bibitem{SY} Y.-T. Siu and S.-T. Yau, {\em Compact K\"ahler manifolds of positive bisectional curvature}, Inven. Math. {\bf 59} (1980), 189-204.

\bibitem{T1} G. Tian, {\em K\"ahler-Einstein metrics with positive scalar curvature}, Invent. Math. {\bf 130} (1997), 1-39.

\end{thebibliography}
\end{document}